# RANDOM WALKS ON A FINITE GROUP AND THE FROBENIUS–SCHUR THEOREM

Olexandr Vyshnevetskiy    Alexander Bendikov

## Abstract

We consider random walk on a finite group $G$ as follows. We can consider $G$ as a group of substitutions. Randomly (i.e. with probability $U(g)=|G|^{-1}$) we choose a substitution $g \in G$ and execute it twice in a row, i.e. execute a substitution $g^2 \in G$. Then the set of squares of elements of the group $G$ be a carrier of a probability $P(g) = \frac{r(g)}{|G|}$ $(g \in G)$, where $r(g)$ is a number of elements $h \in G$ such that $h^2 = g$. Using well-known Frobenius–Schur theorem we find speed of convergence of $n$-fold convolution of $P$ to the uniform probability $U$ and conditions for the convergence.

Let $G$ be a finite group of order $|G|$, $P^{(n)}$ $n$-fold convolution of a probability $P$ on the group $G$. The question of the rate of convergence $P^{(n)}$ to uniform probability on the group $G$ is considered in many papers (see e.g. survey [1] or [2] for groups which are countable sums of finite groups).

For an arbitrary function $f: G \to G$ we denote $M(g) = \{h \in G \mid f(h) = g\}$ and $r(g) = |M(g)|$ $(g \in G)$. If an element $g$ runs through the group $G$, then $f(g)$ takes each value from $f(G) \subset G$ exactly $r(g)$ times. Therefore, for any number function $\chi$ on $G$

$$\sum_g \chi(f(g)) = \sum_t r(t)\chi(t) \qquad (1)$$

(here and below we write $\sum_g$ instead of $\sum_{g \in G}$), where $t$ runs the system of representatives of sets $M(g)$, $g \in G$. Since each element $h \in G$ is included in exactly one of the sets, namely, in the one for which $f(h) = g$, then $\sum_g r(g) = |G|$. Since $r(g) \geq 0$, then

$$P(g) = \frac{1}{|G|} r(g) \qquad (2)$$

is a probability on the group $G$. This is the probability that, at random choice of an element $h \in G$ (i.e., with a uniform (trivial) distribution on $G$), an equality $f(h) = f(g)$ holds for any fixed $g \in G$.

The probability (2) is defined on a group $G$ and its values are non-negative rational numbers with denominators $|G|$. Conversely, any probability $Q$ on $G$ with such values can be given by a formula (2) with an appropriate choice of the function $f$. Namely, let the probability $Q$ takes $m$ non-zero values $Q(g_i) = \frac{n_i}{|G|}$, $i = 1, ..., m$, where $n_i$ are natural numbers, $\sum_{i=1}^{m} n_i = |G|$. Let us divide the group $G$ into $m$ subsets $M_i$, $|M_i| = n_i$, choose an arbitrary element $x_i \in M_i$, $i = 1, ..., m$ and for each $y_i \in M_i$ set $f(y_i) = f(x_i)$.

A function on a group is called a class function if it is constant on each class of conjugate elements of the group. The set $Irr(G) = \{\mathbf{1}_G, \chi_1 ..., \chi_k\}$ of all irreducible complex characters of the group $G$ ($\mathbf{1}_G$ is the principal character) forms a basis of the linear space $L(G)$ of class functions on the group $G$. There is a scalar product in $L(G)$

$$(F_1, F_2) = \frac{1}{|G|} \sum_g F_1(g) \bar{F}_2(g),$$

(where the bar means complex conjugation) and the corresponding norm we denote $\|\cdot\|$. In view of the character orthogonal relations, the set $Irr(G)$ forms an orthonormal basis of the space $L(G)$ with respect to the scalar product.

Any class probability $P \in L(G)$ is expanded on the basis $Irr(G) = \{\mathbf{1}_G, \chi_1 ..., \chi_k\}$, and in view of $\sum_g P(g) = 1$ the coefficient at the main character $\mathbf{1}_G$ is $\frac{1}{|G|}$:

$$P = \frac{1}{|G|} \mathbf{1}_G + \sum_\chi m_\chi \chi \qquad (3)$$

(the last sum is taken over all $\chi \in Irr(G)$, $\chi \neq \mathbf{1}_G$), and the other coefficients

$$m_\chi = (\chi, P) = \frac{1}{|G|} \sum_g \chi(g) P(g)$$

If the probability (2) is a class one, i.e. belongs to $L(G)$, then these coefficients are

$$m_\chi = \frac{1}{|G|^2} \sum_{g \in G} r(g) \chi(g) \qquad (4)$$

We write (3) in the form

$$P = \frac{1}{|G|} \mathbf{1}_G + m_1 \chi_1 + \ldots + m_k \chi_k \qquad (5)$$

where is $k$ the number of non-principal characters of the group $G$.

Let **C**G be a group algebra of a group $G$ over the field **C** of complex numbers. Let us correspond to probability $P(g)$ an element $p = \sum_g P(g) g$ of the algebra **C**G; we denote this element by the same (but small) letter as the function that generats it, and call it *the probability on CG*. Convolution of two functions $P$ and $Q$ on a group $G$

$$(P * Q)(h) = \sum_g P(g) Q(g^{-1} h), \quad h \in G$$

corresponds to the product $pq$ of probabilities on **C**G. Let $U(g) = \frac{1}{|G|}$ be the uniform (trivial) probability on $G$ and $p, u$ the probabilities on the algebra **C**G corresponding to the probabilities $P, U \in L(G)$.

Let $P^{(n)} = P * \ldots * P$ be $n$ - multiple convolution of probability $P$.

Lemma 1 [5].

$$P^{(n)} - U = \frac{1}{|G|} \sum_{j=1}^{k} d_j b_j^n \chi_j \qquad (6)$$

where

$$d_j = \deg \chi_j, \quad b_j = \frac{|G| m_j}{d_j} \quad (j = 1, \ldots k) \qquad (7)$$

Proof. Let $e_j = \frac{d_j}{|G|} \sum_g \chi_j(g) g$ $(j = 1, \ldots k)$. Then from (5) and (7)

$$p = \frac{1}{|G|} \sum_g \mathbf{1}_G g + m_1 \sum_g \chi_1(g) g + \ldots + m_k \sum_g \chi_k(g) g = u + \sum_{j=1}^{k} \frac{|G| m_j}{d_j} e_j = u + \sum_{j=1}^{k} b_j e_j.$$

Since $u, e_1, \ldots, e_k$ are orthogonal idempotents of the center of the algebra

**C**$G$, then
$$p^n = u + \sum_{j=1}^{k} b_j^n e_j,$$

$$p^n - u = \sum_{j=1}^{k} b_j^n e_j = \frac{1}{|G|} \sum_{j=1}^{k} b_j^n \sum_g d_j \chi_j(g) g = \frac{1}{|G|} \sum_g \left( \sum_{j=1}^{k} d_j b_j^n \chi_j(g) \right) g,$$

or, returning to functions on the group $G$, we obtain (6).

Let $G$ be a finite group of permutations of some set, for example, the group of shuffles of a deck of cards. Randomly (i.e. with probability $\frac{1}{|G|}$) we choose a shuffle $g \in G$ and perform it twice in a row, i.e. perform shuffling $g^2 \in G$. The carrier of all such shufflings is the set $T = \{g^2, g \in G\}$. Thus, a random walk is set on the group $G_1 = \langle T \rangle$ generated by $T$, and the corresponding probability is $P(g) = \frac{r(g)}{|G|}$, where $r(g)$ is the number of elements $h \in G$ such that $h^2 = g$ $(g \in G)$.

Taking the function $f(g) = g^2$ in (1) we get
$$\sum_g \chi(g^2) = \sum_t r(t)\chi(t) \qquad (8)$$

Then from (4)
$$|G|m_\chi = \frac{1}{|G|} \sum_g \chi(g^2),$$

A character is called real if all its values are real numbers.

Theorem (Frobenius – Schur) [3, 4]. For any $\chi \in Irr(G)$
$$\frac{1}{|G|} \sum_g \chi(g^2) \in \{0, 1, -1\},$$
where the value $\pm 1$ holds for real characters and only for them.

According to this theorem
$$|G|m_\chi = \begin{cases} 0 \\ \pm 1 \end{cases}, \qquad (9)$$
where the value $\pm 1$ takes place only for real characters.

Let $R(G)$ be the set of all real nonprincipal characters in $Irr(G)$.

Theorem 1. For a random walk defined by probability (2),

$$\|P^{(n)} - U\| = \frac{1}{|G|}\sqrt{\sum_\chi d_\chi^{2(-n+1)}} \tag{10}$$

where the sum is taken over all characters $\chi \in R(G)$.

Proof. In view of (9) and (7) for the considered probability, $b_j = \pm\frac{1}{d_j}$ for characters in $R(G)$ and $b_j = 0$ for the rest. Therefore, (6) can be written as

$$P^{(n)} - U = \frac{1}{|G|}\sum_\chi \pm d_j^{1-n}\chi \tag{11}$$

where the sum is taken over all characters $\chi \in R(G)$. Since the characters in $Irr(G)$ are orthonormal, then applying to (11) Parseval's equality, we obtain (10).

<u>Corollary 1</u>. Asymptotically at $n \to \infty$ $\|P^{(n)} - U\| \sim td^{-n}$, where $d$ is the minimum degree of characters in $R(G)$ and $t$ is the number of such characters of degree $d$.

Let $R_1(G)$ be the set of all linear characters in $R(G)$. Values of characters from $R_1(G)$ are 1 or -1, so its square is the principal character.

<u>Consequence 2</u>. At $n \to \infty$ $\|P^{(n)} - U\| \to 0 \Leftrightarrow R_1(G) = \varnothing$

Recall that $T = \{g^2, g \in G\}$ and $G_1 = \langle T \rangle$. The group $G/G_1$ obviously has exponent 2. Therefore, it is abelian and hence is an elementary abelian 2-group.

Let $X(A)$ be the character group of an abelian group $A$.

<u>Lemma 2</u>. The number $|A|$ is even $\Leftrightarrow R_1(A) \neq \varnothing$.

Proof. The number $|A|$ is even if and only if the group contains a non-trivial involution. Since $A \cong X(A)$, the fact that |A| is even is equivalent to the statement that $X(A)$ containing a non-trivial involution $\alpha$. Character values of $\alpha$ are $\pm 1$ and therefore $\alpha \in R_1(A)$.

<u>Consequence 3</u>. The following conditions are equivalent:
   a) $P^{(n)}$ does not converge to the trivial probability $U$ on the group $G$
   b) $R_1(G) \neq \varnothing$, i.e. $d = 1$
   c) $G \neq G_1$

d) The order of the group $\overline{G} = G/G'$ is even.

Proof. $a) \Leftrightarrow b)$. Consequence 2.

$b) \Leftrightarrow c)$. As was noted above, $G/G_1$ is an Abelian 2-group. $G \neq G_1$ if and only if its order is even. By Lemma 2, this equivalent to $R_1(G/G_1) \neq \varnothing$ $\Leftrightarrow R_1(G) \neq \varnothing$.

$d) \Leftrightarrow b)$. By Lemma 2, the condition $d)$ is equivalent to $R_1(G/G') \neq \varnothing$. Since a character of a factor group is also a character of the group, then $R_1(G) \neq \varnothing$. On the contrary, let $R_1(G) \neq \varnothing$. All the characters in $R_1(G)$ are linear, therefore they contain $G'$ in the kernel and therefore are characters of the factor group $G/G'$. The consequence is proven.

Regardless of whether the conditions of Corollary 3 are satisfied, we have

Corollary 4. The sequence $P^{(n)}$ converges to the trivial probability on the group $G_1$.

Proof. As was proved in [6], for any probability $P$ the sequence $P^{(n)}$ does not converge to $U(\langle V \rangle)$, where $V$ is the support of the probability $P$, if and only if $V$ lies in a nontrivial coset of a group $\langle V \rangle$ by a normal subgroup with a cyclic factor group. However, the carrier $T$ of the probability under consideration contains an element $1 \in G$ and therefore cannot lie in a nonidentity coset. It remains to recall that $G_1 = \langle T \rangle$.

6. Conditions of convergence of a random walk on a finite group. Colloquium Mathematicum , M.S.C., vol. 167, No.1, 2022, p. 109-114.

O. L. VYSHNEVETSKIY
Department of High Mathematics, Kharkiv National Automobile and Highway University,
61002, Ukraine, Kharkiv, st. Yaroslava Mudrogo, 25
*Email address*: alexwish50@gmail.com

A. D. BENDIKOV
Institute of Mathematics, Wroclaw University, Wroclaw, Poland
*E-mail address*: bendikov@math.uni.wroc.pl